\documentclass{amsart}
\usepackage{amsmath,amsthm,amssymb,mathtools,mathrsfs}
\usepackage{geometry}
\usepackage{color}

\newcommand{\A}{{\mathbb A}}
\newcommand{\F}{{\mathbb F}}
\renewcommand{\P}{{\mathbb P}}
\newcommand{\Q}{{\mathbb Q}}
\newcommand{\R}{{\mathbb R}}
\newcommand{\Z}{{\mathbb Z}}
\newcommand{\C}{{\mathbb C}}
\newcommand{\barf}{\bar{f}}
\renewcommand{\bar}[1]{\overline #1}

\newcommand{\cM}{\mathscr{M}}
\newcommand{\cN}{\mathscr{N}}
\newcommand{\cS}{\mathscr{S}}     
\newcommand{\cU}{\mathcal{U}}     
     
\newcommand{\cW}{\mathscr{W}}   

\newcommand{\cH}{\mathscr{H}}

\newcommand{\cNP}{\mathscr{N\!\!P}}

\newcommand{\bG}{\mathbf{G}}

\newcommand{\va}{\vec{a}}
\newcommand{\ba}{\bar{a}}

\newcommand{\denom}{\mathrm{den}}
\newcommand{\heit}{\mathrm{ht}}
\newcommand{\tH}{\tilde{H}}
\newcommand{\tf}{\tilde{f}}
\newcommand{\vn}{\vec{n}}
\newcommand{\vc}{\vec{c}}
\newcommand{\tK}{\tilde{K}}
\newcommand{\chara}{\mathrm{char}}

\newcommand{\hh}{\hspace{1mm}}

\newcommand{\vm}{\overrightarrow{m}}

\newcommand{\sgn}{{\rm{sgn}}}
\newcommand{\NP}{{\rm{NP}}}     
\newcommand{\HP}{{\rm{HP}}}     
\newcommand{\GNP}{{\rm{GNP}}}     
\newcommand{\Tr}{{\rm{Tr}}}

\newcommand{\Gal}{{\rm{Gal}}}

\renewcommand{\hat}{\widehat}
\renewcommand{\tilde}{\widetilde}
\newcommand{\Mat}{\mathrm{Mat}}
\newcommand{\ra}{\rightarrow}
\newcommand{\lra}{\longrightarrow}

\newcommand{\st}{\enskip | \enskip}
\newcommand{\pceil}[1]{\left\lceil #1 \right\rceil}
\newcommand{\pfloor}[1]{\left\lfloor #1 \right\rfloor}

\theoremstyle{plain}
\newtheorem{theorem}{Theorem}[section]
\newtheorem{prop}[theorem]{Proposition}
\newtheorem{lemma}[theorem]{Lemma}
\newtheorem{corollary}[theorem]{Corollary}

\newtheorem{definition}[theorem]{Definition}

\newtheorem*{definition*}{Definition}

\theoremstyle{remark}
\newtheorem{remark}[theorem]{Remark}

\newtheorem*{remark*}{Remark}

\newtheorem{acknowledgments}{Acknowledgments}

\title[Variation of elementary abelian $p$-extensions]
{Asymptotic variation of elementary 
abelian $p$-extensions over $\P^1$}

\author{Hui June Zhu}

\address{
Hui June Zhu,
Department of mathematics,
State University of New York at Buffalo,
Buffalo, NY 14260. 
The United States.
}

\email{hjzhu@math.buffalo.edu}

\date{July 19, 2025}

\keywords{
Elementary abelian $p$-extensions; 
Newton polygons; 
generic Newton polygons; 
Hodge polygons; 
Zeta functions; 
$L$-functions of exponential sums; 
higher rank Artin-Schreier curves.
}
\subjclass{11,14}

\begin{document}

\begin{abstract}
Let $\A^d$ denote the coefficient space 
of all degree-$d$ polynomials $f$ in one variable for some $d\ge 3$.
For any $\bar{f}\in \A^d(\bar\F_p)$,  
a rank-$\ell$ Artin-Schreier curve $X_{\bar{f},\ell}:
y^{p^\ell}-y= \bar{f}$ is called ordinary 
if its normalized Newton polygon 
achieves the infimum in $\A^d(\bar\F_p)$. 
Given $\ell$ and a number field $K$, we show that there exists a Zariski 
dense open subset $\cU$ in $\A^d$, defined over $\Q$, 
such that if $f\in \cU(K)$ then $X_{(f\bmod \wp),\ell}$ is ordinary
for all primes $\wp|p$ with $\deg(\wp)\in\ell\Z$ and $p$ large enough.
\end{abstract}

\maketitle
    
\section{Introduction}

\label{S:intro}

Let $q=p^b$ for some prime $p$ and positive integer $b$. 
A smooth projective curve $C$ over a finite field $\F_q$ has its zeta function 
$Z(C/\F_q,s)$ defined as $Z(C/\F_q,s)=
\exp(\sum_{m=1}^{\infty}\#C(\F_{q^m}) \frac{s^m}{m})$
in variable $s$. It is a rational function (by Weil \cite{We48}),
in fact, a quotient of 
two polynomials in $1+s\Z[s]$.
Its numerator is 
$1+c_1 s+\cdots+c_{2g}s^{2g}=\prod_{i=1}^{2g}(1-\alpha_i s)$,  
where $g$ is the genus of $C$
and reciprocal roots $\alpha_i$ are all Weil $q$-numbers in the sense 
that they are algebraic integers whose complex absolute values of all conjugates 
are equal to $\sqrt{q}$. It then follows that for any prime 
$p'$ different from $p$, the $p'$-adic valuation 
$|\alpha_1|_{p'}=\cdots=|\alpha_{2g}|_{p'} =1$.  
In this paper we study the $p$-adic valuation of these reciprocal roots.

The Newton polygon of $C/\F_q$ is the $q$-adic Newton polygon of
$1+c_1s+\cdots+c_{2g}s^{2g}$, namely, 
the lower convex hull of the points 
$(0,0), \cdots, $ $ (i, v_q(c_i)), \cdots, (2g,v_q(c_{2g}))$ in the real plane $\R^2$,
where $v_q(q)=1$. We denote it by $\NP(C/\F_q)$ or simply $\NP(C)$.
We call this the (normalized) {\em Newton polygon} of $C$.

For each $g\ge 1$, 
consider the set $\cNP_\Z(g)$  
of all Newton polygons $\NP$ as a piece-wise linear function $\NP(x)$
with $x$-coordinate in the interval $[0,2g]$ satisfying the following conditions:
it has endpoints at $(0,0)$ and $(2g,g)$, all vertices $(a,b)$ between the two endpoints have
$a,b\in\Z_{\ge 0}$, and it is symmetric 
in the sense that if there is a line segment of slope $\alpha$ with horizontal length $m_\alpha$ 
then there is a line segment of slope $1-\alpha$ with same horizontal length $m_\alpha$.
We shall also consider the set $\cNP_\Q(g)$ with the same properties except vertices $(a,b)$
have rational coordinates.
In these sets, we order two Newton polygons $\NP\ge \NP'$ if $\NP(x)\ge \NP'(x)$ 
for every $0\le x \le 2g$, we say $\NP$ lies above $\NP'$. 
Note that these sets are partially ordered. Weil theorem says 
that for any projective curve $C/\F_q$ of genus $g$ its Newton polygon lies in $\cNP_\Z(g)$.

Let $\A^d$ denote the coefficient space $\va:=(a_1,\ldots,a_d)$ of 
all polynomials $f=\sum_{i=1}^d a_i x^i$ in one variable $x$.  
For any $\ell\ge 1$, let 
$$
X_{\bar{f},\ell}: y^{p^\ell}-y= \bar{f}
$$ 
be the rank-$\ell$ Artin-Schreier curve 
for some $\bar{f}\in\A^d(\F_q)$.
For rest of the paper, assume that $\ell|b$.
Then  $X_{\bar{f},\ell}$ is an elementary abelian $p$-extension over the projective line $\P^1$.
The natural cover $X_{\bar{f},\ell}\rightarrow \P^1$ defined by sending $(x,y)$ 
to $x$ on points 
is Galois with Galois group $(\Z/p\Z)^\ell$.
In fact, every elementary abelian $p$-extension over 
$\P^1$ that is totally ramified at exactly one point and unramified elsewhere
is of this form (see \cite{GS91}).  
The genus of the curve $X_{\bar{f},\ell}$ is 
$g:=g(X_{\bar{f},\ell})={(p^\ell-1)(d-1)}/{2}$. 

\begin{definition}\label{D:ord}
By Grothendieck's specialization theorem, the following exists
\begin{equation*}
\GNP_{g,\ell,\bar\F_p}:=\inf_{\bar{f}\in\A^d(\bar\F_p)}\NP(X_{\bar{f},\ell}).
\end{equation*}
If $\NP(X_{\bar{f},\ell})=\GNP_{g,\ell,\bar\F_p}$, 
we say $X_{\bar{f},\ell}$ is {\bf ordinary}.
\end{definition}

Our study in this paper was inspired partly by a long-standing conjecture and 
developments for elliptic curves over a number field. 
Given an elliptic curve $E: y^2=x^3+ax+b$ 
with $a,b$ in $\Q$, 
there are infinitely many ordinary primes  by Serre \cite{Se81} 
and there are also infinitely many 
supersingular primes by Elkies \cite{El87}. 
A prime $p$ is {\em ordinary} if 
$\NP(E\bmod p)$ has slopes $0$ and $1$, which is the lower bound.
The set of ordinary primes is of a positive density 
by \cite{Se81} and \cite{El87} independently. 

In this paper, we study a similar question: 
Let $K$ be a number field and let $f\in \A^d(K)$.
We consider whether Artin-Schreier curve $X_{(f\bmod \wp,\ell)}$ is ordinary
at prime $\wp$ of $K$ where $f\bmod \wp$ has meaning,  
that is, $\wp$ is coprime to denominators of coefficients of $f$.
(Note that $f\bmod \wp$ has meaning for almost all $\wp$, 
so we shall tacitly assume so whenever considering the reduction $f\bmod \wp$.) 
We say the prime $\wp$ is {\bf ordinary} for $f$ and $\ell$ if
$X_{(f\bmod \wp),\ell}$ is ordinary.  

For the case $\ell=1$, for any generic $f$,
it is known that $X_{(f\bmod \wp),1}$ is ordinary
when $\wp | p$ for all $p$ large enough (see \cite{Zh03}\cite{Zh04}).
Higher rank case $\ell>1$ is not known to the best of our knowledge, 
and the argument for $\ell=1$ case does not apply to higher rank.

\begin{theorem}
    \label{T:main1}
Let $d\ge 3$.
There exists a Zariski dense open subset $\cU$ of $\A^d$ defined over $\Q$, such that 
if $f\in \cU(K)$ for a number field $K$ then $X_{(f\bmod \wp),\ell}$ is ordinary
for all primes $\wp|p$ of $K$ with $\deg(\wp)\in\ell\Z$ and $p$ large enough.
\end{theorem}

We remark that in the above theorem the lower bound for $p$ depends only on $f$ (independent of $\ell$).

\begin{corollary}\label{C:1}
Let $f\in\cU(K)$ for any number field $K$. 
If $\wp|p$ is a prime of $K$ with $\deg(\wp)\in\ell\Z$ and $p\equiv r\bmod d$ is large enough,
then the slope multi-set of
$\NP(X_{f\bmod \wp,\ell})$ is equal to 
$$
\left\{\frac{1}{d}+\frac{\varepsilon_1}{d(p-1)},\frac{2}{d}+\frac{\varepsilon_2}{d(p-1)},\ldots, \frac{d-1}{d}+\frac{\varepsilon_{d-1}}{d(p-1)}\right\}^{p^\ell-1},$$
where $\varepsilon_i\in \Z$ depends only on $d$ and $r$, and $-(i-1)(d-1)\le \varepsilon_i \le i(d-1)$.
\end{corollary}

\begin{remark}\-
\begin{enumerate} 
\item 
If $f$ is not in $\cU$ or $p$ is small, it is not true that
the Newton polygon of $f$ at each prime $\wp|p$ of degree divisible by $\ell$ 
would achieve the lower generic bound $\GNP_{g,\ell,\bar\F_p}$.
For a counter-example when $p=2$, see Proposition \ref{P:counter-example}
and the discussion in Section \ref{S:2}.
\item 
The theorem shows 
that for any number field $K$ containing a $\Z/\ell\Z$-subextension
a generic $f$ in $\A^d(K)$ has infinitely many ordinary primes. 
See Corollary \ref{C:2}.
\end{enumerate}
\end{remark}

We accomplish the proofs of our main theorems via  Dwork method. 
Applying Dwork methods to this line of study traces back to 
\cite{AS89} and \cite{Wan93}. 
The studies of Newton slopes in Artin-Schreier-Witt towers (see \cite{DWX16} or \cite{KZ18} 
for references) have been extended to the higher-rank case recently 
(see \cite{RWXY} and \cite{Ta22}). 
See also Wan's recent expository article \cite{Wan21} for more related literature and developments.
 
The organization and outline of this paper are as follows:
We introduce some preliminaries of rank-$\ell$ $L$-functions of 
exponential sums in Section \ref{S:2}, and deduce that, as multi-sets of slopes,
$\NP(X_{\bar{f},\ell})$ is the union of Newton polygons of rank-$1$
$L$-functions of $\bar\alpha \bar{f}$
for some $\bar\alpha$ (see \eqref{E:zeta_function}).
We provide an example in which the rank-$1$ 
$L$-functions of $\bar{f}$ and $\bar\alpha \bar{f}$ have distinct Newton polygons for some 
$\bar\alpha$ in characteristic $2$. 
In Section \ref{S:3}, 
we define a {\it deformed Dwork matrix} $\bG$ 
whose `specialization at each point in $\A^d$' yields an original Dwork matrix $G$. 
This allows us to establish 
Proposition \ref{P:L-function} to compute the Newton polygons of $L$-functions.
Separately, we define {\em asymptotic generic Newton polygons}, 
and prepare a sufficient condition for a matrix $G$ to be {\em good enough} 
to achieve them (in Proposition \ref{P:compute!}).
A key new idea of the paper lies in Section \ref{S:4}. Vaguely speaking, 
we need to construct an open set $\cU$ such that if it contains $f$ then it 
contains all its constant multiples
$\alpha f$. This is accomplished in the construction of {\em global generic polynomial} 
for $\A^d$ (see Definition \ref{D:Psi_r}). 
We demonstrate in Proposition \ref{P:est}
how the global generic polynomial (which is independent of $p$) 
allows us to tell which specialization of $\bG$ is {\em good enough} (which is a $p$-adic condition). 
The proof of the main theorems lie in Section \ref{S:proofs}.
We first prove several theorems for $L$-functions
in Theorems \ref{T:sameNP} and \ref{T:main2},
and then derive Theorem \ref{T:main1} easily
after them. In Section \ref{S:6} we study a 1-parameter family of $f=x^d+ax$ and conclude with
Theorems \ref{T:1-parameter} and \ref{T:main1c}, analog to Theorems \ref{T:main2} and \ref{T:main1}, respectively.

\begin{acknowledgments}
The author thanks the referees for their very valuable comments regarding the paper's exposition.
\end{acknowledgments}

\section{Higher rank $L$-function and a counter-example}
\label{S:2}

Let $\bar{f}\in\A^d(\F_q)$, and for the rest of the paper we write
\begin{equation}\label{E:f}
\bar f(x)=\sum_{i=1}^d \bar{a}_i x^i.
\end{equation}
Let $\ell\in\Z_{\ge 1}$ such that $\ell|b$.  
Let $\chi_\ell:\F_{p^\ell}\to \bar\Q_p^*$ be a nontrivial additive character of $\F_{p^\ell}$.
We have $\chi_\ell(1)=\zeta_p$ for some primitive $p$-th root of unity.
We define the $k$-th {\em exponential sum} for $\bar{f}$ over $\F_q$ corresponding to $\chi_\ell$ as \begin{eqnarray}
    S_{\bar f}(k,\chi_\ell) = \sum_{x\in \F_{q^k}} \chi_\ell(\Tr_{\F_{q^k}/\F_{p^\ell}}(\bar f(x))).
\nonumber
\end{eqnarray}
Its $L$-function given by
\begin{eqnarray}
    L_{\bar f}(\chi_\ell,s) = \exp\left(\sum_{k=1}^\infty S_{\bar f}(k,\chi_\ell)\frac{s^k}{k}\right)
\label{E:L-function_1}
\end{eqnarray}
is a polynomial in $1+s\Z_p[\zeta_p][s]$ of degree $d-1$.

\begin{definition}\label{D:NP_q}
The $q$-adic {\em Newton polygon $\NP_q(h(s))$ of 
a power series} $h(s)=1+h_1s+h_2s^2+\ldots$ 
in $1+s\bar\Q_p[[s]]$ is the lower convex hull of points 
$(0,0), (1,v_q(h_1)), (2,v_q(h_2)),\ldots$ in $\R^2$.
For $\bar{f}\in\A^d(\F_q)$, 
we call $\NP_q(L_{\bar{f}}(\chi_\ell,s))$, 
or simply denoted by $\NP(L_{\bar{f}}(\chi_\ell,s))$,
the {\bf (normalized) Newton polygon} of $L_{\bar{f}}(\chi_\ell,s)$.
We write $\NP_q^{<1}$ for its slope $<1$ part.
\end{definition}

\begin{remark}
We  assume $\bar{f}$ in (\ref{E:f}) has no constant term  because 
$L_{\bar f}(\chi_\ell,s)$ and $L_{\bar f+a_0}(\chi_\ell,s)$ have the same Newton polygon. 
In fact, $S_{\bar f+a_0}(k,\chi_\ell) = (\zeta_p^i)^k S_{\bar f}(k,\chi_\ell)$ for some integer $i$ 
that is independent of $a_0$. 
This will simply multiply the reciprocal roots of $L_{\bar f}(\chi_\ell,s)$ by $\zeta_p^i$, 
leaving its $p$-adic Newton polygon unchanged. 
\end{remark}

Let $c_1,\ldots, c_\ell$ be a basis for $\F_{p^\ell}$ over $\F_p$, and let $c_1^*,\ldots, c_\ell^*$ denote its trace dual basis, that is, $\Tr(c_ic_j^*)=\delta_{ij}$.
For $z\in \F_{p^\ell}$ we have $\chi_\ell(z) = \prod_{i=1}^\ell \chi_\ell(c_i^*)^{\Tr_{\F_{p^\ell}/\F_p}(c_i z)}$. Therefore,
\begin{equation}
 \label{eq:S_def}
S_{\bar{f}}(k,\chi_\ell)
 = \sum_{x\in \F_{q^k}} \prod_{i=1}^\ell \chi_\ell(c_i^*)^{\Tr_{\F_{q^k}/\F_p}(c_i \bar f(x))}.
\end{equation}
Then $\chi_\ell(c_i^*) =\zeta_p^{n_i}$ for some $0\leq n_i \leq p-1$.
Write $\alpha \coloneqq \sum_{i=1}^\ell n_ic_i$ in $\F_{p^\ell}$;
Then  (\ref{eq:S_def}) becomes
\begin{equation}\label{E:S_f}
    S_{\bar f}(k,\chi_\ell) = \sum_{x\in \F_{q^k}} \zeta_p^{\Tr_{\F_{q^k}/\F_p}(\alpha \bar f(x))}.
\end{equation}
Combining (\ref{E:L-function_1}) and (\ref{E:S_f}),  we have
$
L_{\bar{f}}(\chi_\ell,s)=L_{\alpha\bar{f}}(\chi_1,s)
$.
This reduces rank-$\ell$ $L$ function to 
rank-$1$ $L$ function. We have proved the following proposition.
We remark here that this manner of reduction was 
well known (for recent other application, see for example 
\cite[Section 3]{Ta22}).

\begin{prop}
\label{P:reduction}
Let $\bar{f}\in\F_q[x]$ and $\chi_\ell: \F_{p^\ell}\lra \bar\Q_p^*$ be a nontrivial character.
Then 
\begin{equation}\label{E:L-L}
L_{\bar{f}}(\chi_\ell,s)=L_{\bar\alpha\bar{f}}(\chi_1,s)
\end{equation}
for some $\bar\alpha\in\F_{p^\ell}^*$ and a nontrivial character $\chi_1:\F_p\lra \bar\Q_p^*$. 
\qed
\end{prop}

Let $Z(X_{\bar{f},\ell}/\F_q,s)$ 
be the Zeta function of the projective 
curve $X_{\bar{f},\ell}$ given by the equation 
$y^{p^\ell}-y=\bar{f}(x)$ over $\F_q$ (that contains $\F_{p^\ell}$); then 
\begin{equation*}
Z(X_{\bar{f},\ell}/\F_q,s)=\frac{\prod_{\bar\alpha\in\F_{p^\ell}^*}L_{\bar\alpha\bar{f}}(\chi_1,s)}{(1-s)(1-qs)}.
\end{equation*}
Considered as slope multi-sets, we hence have
\begin{equation}
 \label{E:zeta_function}
\NP(X_{\bar{f},\ell}/\F_q)=\bigcup_{\bar\alpha\in\F_{p^\ell}^*}\NP_q(L_{\bar\alpha\bar{f}} (\chi_1,s))
\end{equation}
If $\bar\alpha\in\F_p^*$, then 
$L_{\bar\alpha \bar{f}}(\chi_1,s)$ is a conjugate of $L_{\bar{f}}(\chi_1,s)$ 
under the Galois action $\Gal(\Q_p(\zeta_p)/\Q_p)$ (on coefficients only); 
hence, they have the same Newton polygon. 
By \eqref{E:zeta_function}, 
$$Z(X_{\bar{f},1}/\F_q,s)=\frac{\prod_{\sigma\in \Gal(\Q(\zeta_p)/\Q)}\sigma(L_{\bar{f}}(\chi_1,s))}{(1-s)(1-qs)},$$ hence 
$\NP(X_{\bar{f,1}})$ is the dilation of  $\NP(L_{\bar{f}}(\chi_1,s))$ 
by a factor of $p-1$. However, if $\alpha$ lies in $\F_q^*$ but not in $\F_p^*$, 
it is not true
any more that 
$L_{\bar\alpha \bar{f}}(\chi_1,s)$ and $L_{\bar{f}}(\chi_1,s)$ have the same Newton polygon. 
In fact,
we shall provide a counter-example below in which they do not share the same Newton polygon.
(This answers a question of Wan in the one variable case \cite[Conjecture 8.15]{Wan21}. 
See also computation in the case $p=5$ in \cite{Sch23}.)

\begin{prop}
\label{P:counter-example}
Let $a_{11},a_9,a_5\in\bar\F_2^*$ and $a_{11}a_9=a_5^4$.
Let $\bar\alpha\in\bar\F_2^*$ such that $\bar\alpha\ne 1$.
Write $\bar{f}(x)= a_{11}x^{11}+a_{9}x^9+a_5x^5$.
Consider two hyperelliptic curves (we have $\ell=1$ suppressed in notation below for simplicity) 
\begin{eqnarray*}
X_{\bar{f}}:&& y^2-y=\bar{f}(x)\\
X_{\bar\alpha\bar{f}}:&& y^2-y=\bar\alpha \bar{f}(x).
\end{eqnarray*}
Then $\NP(X_{\bar{f}})$ is a straight line of slope $1/2$,
and $\NP(X_{\bar{\alpha}\bar{f}})\neq \NP(X_{\bar{f}})$. 
Moreover,
$$\NP(L_{\bar\alpha\bar{f}}(\chi_1,s))\neq \NP(L_{\bar{f}}(\chi_1,s)).$$
\end{prop}

\begin{proof}
By \cite[Theorem 1]{SZ02}, $X$ is isomorphic to 
$y^2-y=x^{11}+c_3x^3+c_1x$ for some $c_1,c_3\in\bar\F_2$ since $a_{11}a_9=a_5^4$. Hence 
it is supersingular, that is, $\NP(X_{\bar{f}})$ 
is a straight line of slope $1/2$. 
On the other hand, $X_{\bar\alpha\bar{f}}$ can not be supersingular. 
Otherwise, by the same theorem its coefficients satisfy 
$(\bar\alpha a_{11})(\bar\alpha a_9)=(\bar\alpha a_5)^4$. That is ${\bar\alpha}^2 a_{11}a_{9}={\bar\alpha}^4 a_5^4$.
But this implies ${\bar\alpha}^2=1$, which contradicts our hypothesis. Hence $X_{\bar\alpha\bar{f}}$ is not supersingular.  This shows that $\NP(X_{\bar\alpha\bar{f}})\ne \NP(X_{\bar{f}})$.
Therefore, $\NP(L_{\bar\alpha\bar{f}}(\chi_1,s))\neq \NP(L_{\bar{f}}(\chi_1,s))$.
\end{proof}

For primes $p$ that are large, we further investigate the relationship between the two Newton polygons in \S \ref{S:proofs}.

\section{Computational aspect of Newton polygons of $L$-function}
\label{S:3}

Denote by $\Q_q$ the unramified extension of $\Q_p$ of degree $b$, and let $\Z_q$ be its 
ring of integers. We recall $p$-adic Dwork theory and set up basic notations for the rest of the paper. 
Two Newton polygons $\GNP(\A^d,\F_p)$ and $\GNP(\A^d(1),\F_p)$ are explicitly given (in Definition \ref{D:GNP}). We derive a sufficient condition for
the lower bound $\GNP(\A^d,\F_p)$ to be achieved. 
In particular, Proposition  \ref{P:compute!}
will be core ingredients for the proof of Theorem \ref{T:sameNP}.

\subsection{$L$-function via Dwork theory}
\label{S:3.1}

For any complete field $R$ with (non-Archimedean) multiplicative norm $|\cdot|$, 
and $r\in |R|$ and $0<r<1$, define an $R$-algebra
\begin{equation*}
    \cH(R,r) = \left\{\sum_{i=0}^\infty c_i x^i\in R[[x]] \mid 
    \lim_{i\ra \infty}\frac{|c_i|}{r^i} = 0\right\}.
\end{equation*}
This is a Banach $R$-space with respect to the norm 
$|\sum_i c_i x^i|\coloneqq\sup_i(\frac{|c_i|}{r^i})$. 

Let $E(x)$ be the $p$-adic Artin-Hasse exponential function.
 Let $\gamma$ be a root of
$\log E(x)$ in $\bar\Q_p$ with $v_p(\gamma)=\frac{1}{p-1}$ and 
such that $E(\gamma)=\zeta_p=\chi_1(1)$.
Note that $\Z_p[\zeta_p]=\Z_p[\gamma]$ and $\gamma$ is a uniformizer.
Let $\tau \in \Gal(\Q_q(\zeta_p)/\Q_p(\zeta_p))$ be the Frobenius lift.
For any element $\bar{c}\in\F_q$, we write $\hat{c}$ for its Teichm\"uller lifting in $\Z_q$,
note that $\tau(\hat{c})=\hat{c}^p$.
This action extends to an automorphism 
of $\cH(\Q_q(\zeta_p),r)$ by setting $\tau(x)=x$. 

For $\barf(x) = \sum_{i=1}^d \ba_ix^i\in\A^d(\F_q)$, consider the power series in $\Z_q[\zeta_p][[x]]$:
\begin{equation}
    \label{eq:E_f_def}
    E_{\barf}(x) = \prod_{i=1}^{d}E(\hat{a_i}\gamma x^i).
\end{equation}
We shall write $\cH:=\cH(\Q_q(\zeta_p),r)$ for some unspecified $0<r<1$.
Notice $E_{\barf}(x)\in \cH$. 

Write $E(\gamma x)=\sum_{i=0}^\infty\lambda_i x^i$. Then $\lambda_i \in (\Q\cap\Z_p)[\zeta_p]$.
We know $\lambda_i=\frac{\gamma^i}{i!}$ for $1\le i\le p-1$, and $\lambda_0=1$.
Let $A_1,\ldots,A_d$ be variables. Expand $\prod_{i=1}^d E(A_i \gamma x^i)$
as a power series in $(\Z_p[\zeta_p][A_1,\ldots,A_d])[[x]]$ we have
$$\prod_{i=1}^d E(A_i \gamma x^i)=\sum_{n=0}^\infty \bG_n x^n,$$
where $\bG_0=1$ and 
\begin{eqnarray}
\label{E:G_n}
    \bG_n &=& \sum_{\substack{m_i\geq 0\\ \sum_{i=1}^dim_i = n}}\lambda_{m_1}\lambda_{m_2}\cdots \lambda_{m_d}A_1^{m_1}\cdots A_{d}^{m_d} \qquad \mbox{ for $n\ge 1$}.
\end{eqnarray}

\begin{definition}\label{D:deformed}
Assuming $\bG_n=0$ for all $n<0$, we define the {\bf deformed Dwork matrix} for $\A^d$ as
\begin{eqnarray}
\label{E:G}
 \bG &=&(\bG_{pi-j})_{i,j\ge 1}.
\end{eqnarray}
This is a matrix over  $\Z_p[\zeta_p][A_1,\ldots,A_d]$.
Its evaluation at $\hat\va=(\hat{a_1},\ldots,\hat{a_d})$ is 
\begin{eqnarray*}
\bG(\hat\va)&=&(\bG_{pi-j}(\hat\va))_{i,j\ge 1}.
\end{eqnarray*}
\end{definition}

Notice that 
$$E_{\barf} (x)=\prod_{i=1}^{d} E(\gamma \hat{a}_i x^i)
=\prod_{i=1}^d\sum_{k=0}^\infty \lambda_k \hat{a}_i^k x^{ik}
=\sum_{n=0}^{\infty}\bG_n(\hat\va) x^n.
$$
Let $\psi_p(\sum_{i=1}^\infty c_ix^i) = \sum_{i=1}^\infty c_{pi}x^i$.
Define $\psi:\cH\to \cH$ by 
$$
\psi \coloneqq \tau^{-1}\circ \psi_p\circ E_{\barf}(x)
$$ 
where $E_{\barf}(x)$ denotes the multiplication-by-$E_{\barf}(x)$ map.
Notice that  $\psi^b$ is a Dwork operator on the Banach space $\cH$. 
Fix a formal basis $\{1,x,x^2,\ldots\}$ for $\cH$, that is, 
every element of $\cH$ can be uniquely represented as
$\sum_{i=0}^{\infty}c_i x^i $ with $\frac{|c_i|}{r^i}\rightarrow 0$ 
as $i\rightarrow \infty$.
Meanwhile, 
consider the restriction of $\psi_p E_{\bar{f}}(x)$ on subspace $x\cH$ 
generated by $\{x,x^2,\ldots\}$ in $\cH$, this is also its formal basis. 
Write
$$G(\bar{f}):=\Mat(\psi_p E_{\bar{f}}(x)|x\cH).$$
Then it is equal to the evaluation of the deformed Dwork matrix: 
\begin{equation}
\label{E:Gfbar}
G(\bar{f})=(\bG_{pi-j}(\hat\va))_{i,j\ge 1}=\bG(\hat\va).
\end{equation}

To introduce the next result, we shall have a notation in order:
for a square matrix $M$ over $\bar\Q_p$ and an 
automorphism $\tau\in\Gal(\bar\Q_p/\Q_p)$, we denote
(whenever exists) 
\begin{eqnarray}\label{E:define_M}
M_{[b]}=M^{\tau^{b-1}}M^{\tau^{b-2}}\cdots M^\tau M.
\end{eqnarray}

\begin{prop}
\label{P:L-function}
Let $\NP_q$ and $\NP_q^{<1}$ be as in Definition \ref{D:NP_q}.
Let $\tau\in\Gal(\Q_q(\zeta_p)/\Q_p(\zeta_p))$ be the lifting of the Frobenius. 
Then we have
$$\NP_q(L_{\bar{f}}(\chi_1,s))=\NP_q^{<1}(\det(1-G(\bar{f})_{[b]}s))=\NP_q^{<1}(\det(1-\bG(\hat\va)_{[b]}s)).
$$
\end{prop}
\begin{proof}
By \eqref{E:Gfbar}, it remains the show the first equation.
For ease of notation, in this proof (only) we write $G$ for $G(\bar{f})$.
Since $\psi=\tau^{-1}\psi_p E_{\barf}(x)$ we have $\Mat(\psi|x\cH)=G^{\tau^{-1}}$. 
Since $\tau^b=1$, we have 
\begin{equation*}\label{E:psi}
\Mat(\psi^b|x\cH)=G^{\tau^{-1}}G^{\tau^{-2}}\cdots G^{\tau^{-(b-1)}}G^{\tau^{-b}}=G_{[b]}. 
\end{equation*}
Notice $\Mat(\psi)$ has the first row equal to $(1,0,0,\ldots)$, so does $\Mat(\psi^b)$.
By the equation above, it follows that 
\begin{equation}\label{E:det_22}
\det(1-s\psi^b)=(1-s)\det(1-s\psi^b|x\cH)
=(1-s)\det(1-G_{[b]}s).
\end{equation} 
It is proved by Dwork (\cite{Dw60}) that the determinant is well-defined 
and independent of choice of matrix representations of the Dwork operator (i.e.,
independent of the choice of formal basis).  
Note that Dwork operator is a completely continuous operator in \cite{Se62}. 
On the other hand, by Dwork theory we have
\begin{eqnarray}\label{E:Dwork}
L_{\barf}(\chi_1,s) & = &
\frac{\det(1-s \psi^b)}{(1-s)\det(1-s\psi^b)}.
\end{eqnarray}
By \eqref{E:det_22} and \eqref{E:Dwork}, we have
\begin{align*}
L_{\barf}(\chi_1,s)=\frac{(1-s)\det(1-G_{[b]}s)}{(1-s)(1-qs)\det(1-G_{[b]}qs)}
=
\frac{\det(1-G_{[b]}s)}{(1-qs)\det(1-G_{[b]}qs)}.
\end{align*}
Thus we get the following factorization in
$1+s\Z_q[\zeta_p][[s]]$:
$$
\det(1-G_{[b]}s)=L_{\barf}(\chi_1,s)\cdot (1-qs)\cdot \det(1-G_{[b]}qs).
$$
Notice that $q$-adic slopes of factors $(1-qs)$ and $\det(1-G_{[b]}qs)$ 
are all $\ge 1$. But it is well known that 
$L_{\bar{f}}$ is a polynomial of degree $d-1$ (see \cite{Bom66} or \cite{AS89}), 
whose Newton polygon slopes lie in the open interval $(0,1)$ with endpoints $(0,0)$
and $(d-1,(d-1)/2)$. Thus $\NP_q(\det(1-G_{[b]}s))$ has a vertex at $(d-1,(d-1)/2)$,
and  
$\NP_q(L_{\bar{f}}(\chi_1,s))=\NP_q^{<1}(\det(1-G_{[b]}s))$.
\end{proof}

\subsection{Asymptotic generic Newton polygon}

For any matrix $M=(m_{ij})_{i,j\ge 1}$ over $\bar\Q_p$, 
define a degree $d-1$ auxiliary polynomial of $M$ in $1+s\bar\Q_p[s]$ as follows  
\begin{equation}
\label{E:P_Ms}
Q_M(s)\coloneqq 1+\sum_{n=1}^{d-1} \det((m_{ij})_{1\le i,j\le n})s^n.
\end{equation}
We say $M=(m_{ij})_{i,j\ge 1}$ 
is {\it nuclear} if $\lim\limits_{i\ra \infty}\inf_j v_p(m_{ij}) = \infty.$
(This definition is the same as Serre's \cite{Se62}, 
we opt for the concrete definition for simplicity.)

\begin{prop}
\label{P:transform}
Let $\delta\in\R_{>0}$ be fixed. Let $t\ge 2$.
Let $h_1< h_2<\ldots <h_t<1\le h_{t+1}\le h_{t+2}\le \cdots 
$ be a real number sequence such that 
$h_{i+1}-h_i\ge \delta$ for every $1\le i\le t$.
Suppose $M=(m_{ij})_{i,j\ge 1}$ is any nuclear matrix over $\bar\Q_p$ satisfying the following properties:
\begin{enumerate}
\item $\inf_{j\ge 1}v_p(m_{ij})\ge h_i$ for every $i\ge 1$, 
\item $v_p(\det (m_{ij})_{1\le i,j\le n})< \sum_{i=1}^{n}h_i + \frac{\delta}{2}$ for every $1\le n\le t$.
\end{enumerate}
Then 
$$
\NP^{<1}_q(\det(1-M_{[b]}s))=\NP^{<1}_p(\det(1-Ms))=\NP_p(Q_M).
$$

\end{prop}

\begin{proof}
We first prove the second equality.
Writing $\det(1-Ms)=\sum_{k=0}^\infty c_k s^k$, we have 
$(-1)^n c_n=\sum_{|\cS|=n} \det((m_{i,j})_{i,j\in\cS})$,
where $\cS$ ranges over all sets of $n$ positive integers $k_1<k_2<\cdots<k_n$.
It suffices to consider only slope $<1$ part of the Newton polygon so by hypothesis (1), 
we may assume $1\le n\le t$ for the rest of this argument.
If $\cS\ne \{1,2,\ldots,n\}$, then 
$v_p(\det( m_{ij})_{i,j\in\cS})\ge h_{k_1}+h_{k_2}+\cdots +h_{k_n}
\ge h_1+h_2+\cdots + h_{n-1}+h_{n+1}= \sum_{i=1}^n h_i+(h_{n+1}-h_n)
>v_p(\det (m_{ij})_{1\le i,j\le n})$ by the hypotheses.
Hence 
$v_p(c_n)=v_p \det((m_{ij})_{1\le i,j\le n})$. Our second equality follows.

Write
$\det(1-M_{[b]}s)=\sum_{k=0}^{\infty}C_k s^k$. 
Our hypotheses (1) and (2) imply that $$v_p(\det(m_{ij})_{1\le i,j\le n})
<\sum_{i=1}^n h_i +\frac\delta{2}\le \sum_{i=1}^n h_i+\frac{h_{n+1}-h_n}{2}.$$
Thus it follows from \cite[Theorem 5.5]{LZ04}
that $v_q(C_n)=v_p(\det (m_{ij})_{1\le i,j\le n})$. This proves the first equality.
\end{proof}

In this paper let 
\begin{eqnarray}
(\Z/d\Z)^*:=\{c\in\Z\mid 0< c\le d-1, (c,d)=1\}.
\end{eqnarray}
 We always write $(c\bmod d)$ for
the representative of $c$ in $\{0,1,\ldots,d-1\}$. 
Fix $r\in(\Z/d\Z)^*$. Set
\begin{eqnarray}
r_{ij}&\coloneqq&(-(ri-j)\bmod d),
\label{E:r_ij}
\\
M_n&\coloneqq&\min_{\sigma\in S_n}\sum_{i=1}^n r_{i,\sigma(i)}.
\label{E:M_n}
\end{eqnarray}

For convenience of the reader, we recall below 
two Newton polygons $\GNP(\A^d,\F_p)$ and $\GNP(\A^d(1),\F_p)$, first introduced in 
\cite{Zh03}. 
(Note that $\A^d(1)$ denotes the subfamily of all polynomials $f=x^d+ax$, which we shall study in Section \ref{S:6}.) 

\begin{definition}
\label{D:GNP}
Let $r\in (\Z/d\Z)^*$ and $p$ be any prime such that 
$p\equiv r\bmod d$. \\
(1) Let $\GNP(\A^d,\F_p)$ be the Newton polygon with endpoints and vertices at 
$$(0,0), (1,w_1), (2,w_2), \ldots, (d-2,w_{d-2}), (d-1,\frac{d-1}{2}),$$ where
\begin{equation}\label{E:GNP}
w_n\coloneqq\frac{n(n+1)}{2d}+ \frac{M_n}{d(p-1)}.
\end{equation}
We call this the {\bf asymptotic generic Newton polygon 
for $\A^d$}. 
Notice that for $p\rightarrow \infty$ we have $w_n\rightarrow \frac{n(n+1)}{2d}$, 
hence these are indeed vertices for large $p$ (which will be the case in this paper for our purpose).

(2) 
Let $r'_{ij}=(ri-j\bmod d)$ and let $M_n'=\min_{\sigma\in S_n}\sum_{i=1}^n r_{i,\sigma(i)}'$.
Let $\GNP(\A^d(1), \F_p)$ be the Newton polygon with vertices
at 
$$
(0,0), (1,w_1'), (2,w_2'), \cdots, (d-2,w_{d-2}'), (d-1, \frac{d-1}{2}),
$$
where 
$$
w_n':= \frac{n(n+1)}{2d} + \frac{(d-1)M_n'}{d(p-1)}.
$$
We call this the {\bf asymptotic generic Newton polygon for $\A^d(1)$}. 
\end{definition}

\begin{prop}
\label{P:compute!}
Let $p$ be a prime such that $p>2d^2-2d+1$.
Write $q=p^b$ for some $b\ge 1$. Suppose a  matrix 
$G\coloneqq(G_{ij})_{i,j\ge 1}$ with $G_{ij}\in\Q_q(\zeta_p)$ satisfies  
\begin{enumerate}
\item $v_p(G_{ij})\ge 
\frac{pi-j}{d(p-1)}$ for all $i,j\ge 1$,
\item $v_p(G_{ij})=\frac{\pceil{\frac{pi-j}{d}}}{p-1}$ 
(resp. $\frac{\pfloor{\frac{pi-j}{d}}+r'_{ij}}{p-1}$) for all $i,j\le d-1$,  
\item $v_p(\det (G_{ij})_{1\le i,j\le n})=w_n$ (resp. $w_n'$) 
for all $1\le n\le d-1$.
\end{enumerate}
Then 
$$
\NP_q^{<1}(\det(1-G_{[b]}s))
=\GNP(\A^d,\F_p) 
\mbox{\quad (resp. $\GNP(\A^d(1),\F_p)$)}.
$$
\end{prop}

\begin{proof}
We shall prove the first case for $\GNP(\A^d,\F_p)$, 
and omit the proof of the second for $\GNP(\A^d(1),\F_p)$ 
as they are identical in argument.

Note that $\gamma\in\Q_p(\zeta_p)$ is the root of $\log E(x)$ with $v_p(\gamma)=1/(p-1)$.
Pick and fix a $\gamma^{\frac{1}{d}}\in \bar\Q_p$.
Let 
$M=(m_{ij})_{i,j\ge 1}$ where 
$m_{ij}=G_{ij}\gamma^{\frac{j-i}{d}}$. It is clear that 
\begin{equation}\label{E:MG}
v_p(\det ((m_{ij})_{1\le i,j\le n}))=v_p(\det ((G_{ij})_{1\le i,j\le n}))=w_n.
\end{equation}
By definition of \eqref{E:P_Ms}, we have  
\begin{equation}\label{E:P_M}
\NP_p(Q_M) = \GNP(\A^d,\F_p).
\end{equation}

Write $h_i=\frac{i}{d}$ for 
all $i\ge 1$ and let $\delta=\frac{1}{d}$. 
We shall show that $M$ satisfies the two hypotheses of Proposition \ref{P:transform}.
By our hypotheses (1) and (2), 
\begin{equation}
\label{E:v_p}
v_p(m_{ij})\ge \frac{pi-j}{d(p-1)}+\frac{j-i}{d(p-1)}= h_i \mbox{ for all $i,j\ge 1$}.
\end{equation} 
Hence $M$ is nuclear.
Since $n\le d-1$, $M_n\le n(d-1)$, and $p>2d^2-d+1$, we have
$$w_n=\frac{n(n+1)}{2d}+\frac{M_{n}}{(p-1)d}\le\frac{n(n+1)}{2d}+\frac{(d-1)n}{(p-1)d}
<\sum_{i=1}^n h_i +\frac{d-1}{p-1}<\sum_{i=1}^n h_i +\frac{1}{2}\delta.$$ 
By our hypothesis (3), we have 
\begin{eqnarray}
\label{E:value}
v_p(\det ((G_{ij})_{1\le i,j\le n})=w_n
< \sum_{i=1}^n h_i +\frac{1}{2}\delta.
\end{eqnarray}
By  \eqref{E:MG} and \eqref{E:value}, 
we have
\begin{equation}\label{E:v_p2}
v_p(\det ((m_{ij})_{1\le i,j\le n})<\sum_{i=1}^n h_i +\frac{1}{2}\delta.
\end{equation}
In \eqref{E:v_p} and \eqref{E:v_p2}, we have shown that the hypotheses of Proposition \ref{P:transform} are 
satisfied, thus 
\begin{equation}\label{E:22}
\NP_q^{<1}(\det(1-M_{[b]}s))
=\NP_p(Q_M).
\end{equation}
On the other hand, Fredholm theorem gives that 
$\det(1-G_{[b]}s)=\det(1-M_{[b]}s)$, hence
\begin{equation}\label{E:3}
\NP_q^{<1}\det(1-G_{[b]}s)= \NP_q^{<1}(\det(1-M_{[b]}s)).
\end{equation}
Combining \eqref{E:P_M}, \eqref{E:22} and \eqref{E:3}, we obtain our desired equation
$
\NP_q^{<1}(1-G_{[b]}s)=\GNP(\A^d,\F_p)
$
\end{proof}

\section{Global generic polynomials that are independent of $p$}

\label{S:4}

Let $\bar\Z_p$ denote the integral closure in $\bar\Q_p$.
Henceforth we tacitly pick an embedding from $\bar\Q$ to $\bar\Q_p$ and $\C_p$
once and for all. 

This section is the dry technical part of the paper, 
which the reader may opt to skip at first reading. 
In \S \ref{S:4.1} we construct a nonzero 
polynomial $\tilde\Psi_r$ in $\Q[A_1,\ldots, A_d]$ for each $r$
in $(\Z/d\Z)^*$, 
which is used to define a Zariski dense open subset in $\A^d$.  
In \S \ref{S:4.2}, we prove Proposition \ref{P:est} that shows how 
$\tilde\Psi_r$ controls the $p$-adic valuations $\bG_{ij}(\vc)$ and $\det(\bG_{ij}(\vc))_{i,j\le n}$
if $p\equiv r\bmod d$ for some $\vc\in{\bar\Z_p}^d$. We assume $1\leq i,j,n\le d-1$ below.
 
\subsection{Global generic polynomial}
\label{S:4.1}

We fix $r\in (\Z/d\Z)^*$ below unless otherwise declared.
First recall $r_{ij}=(-(ri-j) \mod d)$ from \eqref{E:r_ij} 
and let $r_{i1}'\coloneqq(ri-1 \mod d)$.
Write $\vm=(m_1,\ldots,m_d)$ below. 
Consider the following nonempty set
\begin{equation*}
    \cM_{ij} \coloneqq \{\vm\in\Z_{\ge 0}^d\st  \sum_{k=1}^{d-1} 
    (d-k) m_k = r_{ij}, \sum_{k=1}^d m_k=d\}.
\end{equation*}
Let 
\begin{equation*}
    \delta_{ij} \coloneqq \begin{cases}
        0 & \mbox{if $j< r_{i1}'+1$}\\
        1 & \mbox{if $j\geq r_{i1}'+1$}.
    \end{cases}
\end{equation*}
For any $\vm\in \cM_{ij}$
we define a nonzero rational number
\begin{equation*}
    h_{\vm,i,j} = \frac{\prod_{\ell=0}^{\sum_{k=1}^{d-1}m_k +\delta_{ij}-1}(\frac{r_{i1}-1}{d}-\ell)}{\prod_{k=1}^{d-1}m_k!}.
    \end{equation*}
Let $A_1,\ldots,A_{d}$ be variables as before.
Let   
 \begin{equation}
 \label{E:H}
\tH_{ij} \coloneqq  \sum_{\vm\in\cM_{ij}} h_{\vm,i,j} A_1^{m_1}\cdots A_d^{m_d}.
\end{equation}
  
\begin{lemma}\label{L:H}
Then $\tilde{H}_{ij}$ is nonzero homogenous of total degree $d$ in $\Q[A_1,\ldots,A_d]$
and is supported on each $\vm\in\cM_{ij}$, that is, $h_{\vm,i,j}\in\Q^*$ for every $\vm\in\cM_{ij}$.
\end{lemma}
\begin{proof}
Since $d>r_{ij}$ the set $\cM_{ij}$ is not empty.
By hypothesis, $(r,d)=1$ and $i\le d-1$, we have 
$r_{i1}-1\equiv -(ri-1)-1=-ri \not\equiv 0\bmod d$ 
hence $\frac{ri-1}{d}\not\in \Z$. 
Thus $h_{\vm,i,j}\ne 0$ for every $\vm\in\cM_{ij}$.
By construction, $\tH_{ij}$ is homogenous of total degree $d$ in 
$\Q[A_1,\ldots,A_d]$.
\end{proof}
 
Consider the determinant 
\begin{equation}\label{E:det2}
\det (\tH_{ij})_{1\le i,j\le n}=\sum_{\sigma\in S_n}\sgn(\sigma)\prod_{i=1}^n \tH_{i,\sigma(i)}.
\end{equation}
 Let  $S_n^{\min}\coloneqq
    \left\{ \sigma \in S_n | 
    \sum_{i=1}^n r_{i,\sigma(i)} = M_n \right\}$, where $M_n$ is given in \eqref{E:M_n}.
Let 
\begin{eqnarray}\label{E:f_n}
   \tf_n & \coloneqq&\sum_{\sigma \in S_n^{\min}}\sgn(\sigma)\prod_{i=1}^n \tH_{i,\sigma(i)}.
\end{eqnarray}
This is clearly a polynomial in $\Q[A_1,\ldots,A_d]$, as each summand is due to Lemma \ref{L:H}.
But it is highly non-trivial to tell whether it is nonzero, which we shall prove in the lemma below.
Some comments of our strategy are in order: 
We consider the polynomials \eqref{E:det2} and \eqref{E:f_n} as sums over $\sigma\in S_n$
and $\sigma\in S_n^{\min}$, respectively. Each $\sigma$-summand 
$\sgn(\sigma)\prod_{i=1}^n\tH_{i,\sigma(i)}$ is a nonzero homogenous polynomial of total degree $nd$ by Lemma \ref{L:H}, the only problem remains: their sum has potential cancellations. 
So it suffices to establish the existence of a monomial term in a $\sigma$-summand of $\tf_n$
that is of a unique monomial order (below we use the lexicography order), 
since its uniqueness eliminates any possibility of cancellation.

 \begin{lemma}
\label{L:f_n} 
The polynomial $\det(\tH_{ij})_{1\le i,j\le n}\in\Q[A_1,\ldots,A_d]$ in \eqref{E:det2} has a unique 
monomial term $\Xi$ of the highest lexicographic order with $A_d>\ldots >A_1$,
 and $\Xi$ lies in $\tf_n$. In particular, $\tf_n$ is a nonzero homogenous polynomial in $\Q[A_1,\ldots,A_d]$ of total degree $nd$.
 \end{lemma}

 \begin{proof}
(1) 
Let $\Delta_{ij}\coloneqq A_d^{d-1}A_{d-r_{ij}}$. Then 
\begin{equation}\label{E:det1}
\det(\Delta_{ij})_{1\le i,j\le n}= A_d^{(d-1)n}\sum_{\sigma\in S_n}\sgn(\sigma)
\prod_{k=0}^{d-1}A_{d-k}^{\#\{1\le i\le n|r_{i,\sigma(i)}=k\}}.
\end{equation}
We first define $\sigma_0\in S_n$ and show that 
the $\sigma_0$-term in the expansion \eqref{E:det1} 
is of the highest lexicographic order, and hence is a unique term.  
We define $\sigma_0\in S_n$ as follows:
Let  $\cS_0=\{(i,j)\in\{1,\ldots, n\}^2\st r_{ij}=0\}$. 
Then inductively, for every $k\ge 1$ let 
$\cS_k=\{(i,j)\in \{1,2\ldots,n\}^2 \st (i,\star), (\star,j)\not\in \cup_{\ell=0}^{k-1} \cS_\ell, r_{ij}=k\}$. 
Notice that $\cS\coloneqq\cup_{k=0}^{d-1}\cS_k$ consists exactly $n$ pairs. Define 
$\sigma_0(i)=j$ for every $(i,j)\in \cS$. 
Since $(r,d)=1$, the map $\sigma_0\in S_n$.
Notice $\# \{1\le i\le n|r_{i,\sigma(i)}=0\}\le \# \cS_0$
and the equality holds precisely when $\sigma(i)=\sigma_0(i)=j$ for all $(i,j)\in\cS_0$. 
In this case, $\{1\le i\le n|r_{i,\sigma(i)}=1\}\le \# \cS_1$
and equality holds precisely when $\sigma(i)=\sigma_0(0)=j$ for all $(i,j)\in \cS_1$, and inductively for all $1\le k\le n$. 
This proves that the highest lexicographic order is achieved precisely at the 
$\sigma_0$-term in \eqref{E:det1}, 
and its monomial is equal to $\Xi\coloneqq A_d^{(d-1)n+\#\cS_0} A_{d-1}^{\#\cS_1}\cdots A_1^{\#\cS_{d-1}}$.
 
(2) 
Let $\eta=\prod_{i=1}^n \eta_{\sigma,i}$ be a  highest lexicographic order monomial in the formal expansion \eqref{E:det2}, where $\eta_{\sigma,i}$ is a monomial in $\tH_{i,\sigma(i)}$.
Then $\eta_{\sigma,i}$ must be of the highest lexicographic order in $\tH_{i,\sigma(i)}$ for each $1\le i\le n$. 
By part (1), $\eta_{\sigma,i}=\Delta_{i,\sigma(i)}$. Combined with part (1), we conclude that 
$\sigma=\sigma_0$. Thus $\Xi$ is the unique monomial of the highest lexicographic order 
in $\det(\tH_{ij})_{1\le i,j\le n}$.

(3) We claim that $\Xi$ lies in $\tf_n$, that is, $\sigma_0\in S_n^{\min}$. 
By our construction of $\sigma_0$, for every $i$ with $r_{ij}=0$ for some $1\le j\le n$, 
 we have $r_{i,j}-r_{i,\sigma_0(i)}=r_{i,j}\ge j-\sigma_0(i)$; 
 For every $i$ with $r_{ij}\ne 0$ for all $1\le j\le n$
 then $r_{i,j}$ is strictly increasing and $r_{i,j}-r_{i,\sigma_0(i)}= j-\sigma_0(i)$ clearly.
Sum over all $i$ we have for all $j\in\{1,\ldots,n\}$
$\sum_{i=1}^{n}(r_{i,j}-r_{i,\sigma(i)})\ge  \sum_{i=1}^{n}(j-\sigma_0(i))$.
For any $\sigma\in S_n$, substitute $j=\sigma(i)$ to the above, we have
$\sum_{i=1}^n r_{i,\sigma(i)} - \sum_{i=1}^{n}r_{i,\sigma_0(i)}\ge 0$.
This proves that $\sigma_0\in S_n^{\min}$.

Since  $\Xi$ is a term that can not be cancelled with other terms in $\tf_n$, we
conclude that $\tf_n$ is a nonzero polynomial in $\Q[A_1,\ldots,A_d]$. It is clearly homogenous of total degree $nd$.
 \end{proof}

\begin{definition}\label{D:Psi_r}
For each $r\in(\Z/d\Z)^*$, let 
$\tH_{r,i,j}$ and $\tf_{r,n}$ denote $\tH_{ij}$ and $\tf_n$ 
as given in \eqref{E:H} and \eqref{E:f_n}, respectively.
We call the following polynomial in $\Q[A_1,\ldots,A_d]$ 
\begin{equation*}
\widetilde{\Psi}_r\coloneqq A_d\prod_{i=1}^{d-1}\prod_{j=1}^{d-1}\tilde{H}_{r,i,j}\prod_{n=1}^{d-2} \widetilde{f}_{r,n}
\end{equation*}
the {\bf global generic polynomial for $\A^d$}. 
We call each factor on the right-hand-side, that is, 
$A_d, \tH_{r,i,j}, \tf_{r,n}$, a {\bf factor} of $\tilde{\Psi}_r$.
\end{definition}

\begin{prop}
\label{P:open} \-
\begin{enumerate}
\item 
Let $\tilde\cW_r$ be the complement
of $V(\tilde{\Psi}_r)$ in $\A^d$. 
Then $\tilde\cW_r$ is a Zariski dense open subset in  
$\A^d$ defined over $\Q$, depending only on $d$ and $r$. 
\item For $p\equiv r\bmod d$ and $p>d$, each factor of $\tilde\Psi_r$ 
lies in $(\Q\cap\Z_p)[A_1,\ldots,A_d]$. 
\end{enumerate}
\end{prop}

\begin{proof}
(1) By Lemmas \ref{L:H} and \ref{L:f_n}, $\tilde{\Psi}_r$ is a non-zero polynomial 
in $\Q[A_1,\ldots,A_d]$. 
Then
$\tilde\cW_r$ is open nonempty and hence dense because  $\A^d$ is irreducible. 
The dependency follows from the construction of $\tilde\Psi_r$.

(2) Notice that for $p>d$ and $p\equiv r\bmod d$, $p$ is coprime to denominators of all coefficients of $\tH_{r,i,j}$ and hence also to those of $\tf_{r,n}$. Hence coefficients of $\tilde{\Psi}_r$ are in $\Z_p$. 
\end{proof}

\subsection{Local generic polynomial at each $p$}
\label{S:4.2}

For this subsection, $p$ is  a prime number  
such that $p\equiv r\bmod d$.

Let 
\begin{equation}
\label{E:Kij}
   \tK_{ij} = \sum_{\vm\in\cN_{ij}} \frac{A_1^{m_1}\cdot\cdots\cdot A_{d}^{m_{d}}}{m_1!\cdots m_d!}
\end{equation}
where 
\begin{equation}\label{E:tM}
\cN_{ij}= \{\vm\in\Z_{\ge 0}^d
\st 
\sum_{k=1}^{d} km_k=pi-j, \sum_{k=1}^{d} m_k=\pceil{\frac{pi-j}{d}}\}.
\end{equation}
Notice that $\tK_{ij} $ is homogenous of total degree $\pceil{\frac{pi-j}{d}}$ by its definition. 
Consider the monomial terms in the expansion
$\det (\tK_{ij})_{1\le i,j\le n}=\sum_{\sigma\in S_n}\sgn(\sigma)\prod_{i=1}^n \tK_{i,\sigma(i)}$.
For $p$ large enough notice that $\tK_{ij}\in\Z_p[A_1,\ldots,A_d]$.

\begin{lemma}\label{L:key2}
Suppose $p\ge (d^2+1)(d-1)$. 
Then $\tK_{ij}\in (\Q\cap \Z_p)[A_1,\ldots,A_d]$ and 
$$
\tK_{ij}\equiv v_i A_d^{\pceil{\frac{pi-j}{d}}-d}\tH_{ij}\pmod p
$$
where $v_i\in \Q\cap\Z_p^*$ depends only on $i$.
\end{lemma}

\begin{proof}
Notice that $\pceil{\frac{pi-1}{d}}\equiv \frac{r_{i1}-1}{d}\pmod p$
and for all $j\ge 1$ we have $\pceil{\frac{pi-j}{d}}\equiv \frac{r_{i1}-1}{d}-\delta_{ij}
\pmod p$. Therefore, 

\begin{align*}
h_{\vm,i,j}
&\equiv \frac{\pceil{\frac{pi-1}{d}}(\pceil{\frac{pi-1}{d}}-1)
\cdots (\pceil{\frac{pi-j}{d}}-(\sum_{k=1}^{d-1}m_k-1))}{m_1!\cdots m_{d-1}!}
\\
&\equiv \frac{\pceil{\frac{pi-1}{d}}!}{m_1!\cdots m_{d-1}!(\pceil{\frac{pi-j}{d}}-\sum_{k=1}^{d-1}m_k )!} \pmod p.
\end{align*}
Then we have 
\begin{eqnarray}
\label{E:HH}
\tH_{ij} 
&=& \sum_{\vm\in\cM_{ij}} h_{\vm,i,j} A_1^{m_1}\cdots A_d^{m_d}\\
&\equiv &
\sum_{\vm\in\cM_{ij}} \frac{\pceil{\frac{pi-1}{d}}!A_1^{m_1}\cdots A_d^{d-\sum_{k=1}^{d-1}m_k}}{m_1!\ldots m_{d-1}!(\pceil{\frac{pi-j}{d}}-\sum_{k=1}^{d-1}m_k)!}
\pmod p.
\nonumber
\end{eqnarray}
Notice 
for  every $\vm=(m_1,\ldots,m_{d-1},m_d)\in\cM_{ij}$
we have $(m_1,\ldots,m_{d-1}, \pceil{\frac{pi-j}{d}}-\sum_{k=1}^{d-1}m_k)\in \cN_{ij}$. 
In fact, this map gives a bijection of the sets $\cM_{ij}$ and $\cN_{ij}$.
Write $v_i=\frac{1}{\pceil{\frac{pi-1}{d}}!}$. It is clearly a rational number
and is a $p$-adic unit.
Then we have 
\begin{eqnarray}
\label{E:KK}
\tK_{ij} &= & 
\sum_{\vn\in \cN_{ij}}\frac{A_1^{n_1}\cdots A_d^{n_d}}{n_1!\ldots n_d!}\\
\nonumber 
&\equiv &
\sum_{\vm\in \cM_{ij}}
\frac{
A_1^{m_1}\cdots A_{d-1}^{m_{d-1}}A_d^{(\pceil{\frac{pi-j}{d}}-\sum_{k=1}^{d-1}m_k)}
}
{m_1!\ldots m_{d-1}!(\pceil{\frac{pi-j}{d}}-\sum_{k=1}^{d-1}m_k)!}
\\
\nonumber &\equiv &
\frac{1}{\pceil{\frac{pi-1}{d}} !}A_d^{\pceil{\frac{pi-j}{d}}-d} 
\sum_{\vm\in \cM_{ij}}
\frac{\pceil{\frac{pi-1}{d}} !       A_1^{m_1}\cdots A_{d-1}^{m_{d-1}}A_d^{d-\sum_{k=1}^{d-1}m_k}
}
{m_1!\ldots m_{d-1}!(\pceil{\frac{pi-j}{d}}-\sum_{k=1}^{d-1}m_k)!}
\\
\nonumber
&\equiv &
v_i A_d^{\pceil{\frac{pi-j}{d}}-d}\tH_{ij} \pmod p,
\end{eqnarray}
where the last congruence follows from \eqref{E:HH}.
\end{proof}

Recall $\bG_n$ is a polynomial in $(\Z_p[\zeta_p])[A_1,\ldots,A_d]$ 
from \eqref{E:G_n}, and $\gamma$ is the uniformizer of $\Z_p[\zeta_p]$ introduced in Section 
\ref{S:3.1}.
Let  $w_n$ be as in \eqref{E:GNP}. 

\begin{lemma}\label{L:congruence}
Let $p>(d^2+1)(d-1)$. Then we have  
$$
\left\{
\begin{array}{llll}
 \bG_{pi-j}&\equiv & 
 v_i\gamma^{\pceil{\frac{pi-j}{d}}}
 A_d^{\pceil{\frac{pi-j}{d}} -d}\tilde{H}_{ij} 
&\pmod{\gamma^{\pceil{\frac{pi-j}{d}}+1}}, \\
 \det ((\bG_{pi-j})_{1\le i,j\le n}) 
 &\equiv  &
 u_n\gamma^{(p-1)w_n}A_d^{N_n}\widetilde{f}_n 
 &\pmod {\gamma^{(p-1)w_n+1}},
    \end{array}
\right.
$$
for some $u_n\in\Z_p^*$ and 
$N_n=(p-1)w_n-nd$.   
\end{lemma}

\begin{proof}   
In the definition \eqref{E:G_n} of $\bG_{pi-j}$,
since $m_1+\cdots +dm_d=pi-j$ and all $m_i\ge 0$, we have 
$\sum_{i=1}^d m_i \ge \pceil{\frac{pi-j}{d}}$. If the equality holds 
then $0\le m_i\le p-1$, and hence 
 we have $\lambda_{m_i}=\frac{\gamma^{m_i}}{m_i!}$ for all $i$.
Thus,
$$\bG_{pi-j} \equiv \gamma^{\pceil{\frac{pi-j}{d}}}\sum_{\vm\in \cN_{ij}}
\frac{\prod_{k=1}^{d} A_k^{m_k}}{m_1!\cdots m_d!} 
\equiv \gamma^{\pceil{\frac{pi-j}{d}}}\tK_{ij} \pmod 
{\gamma^{\pceil{\frac{pi-j}{d}}+1}}.
        $$
Combined with Lemma \ref{L:key2}, $v_p(p)=v_p(\gamma^{p-1})$,
we obtain the first congruence.

To prove the second congruence, first note that
\begin{equation*}
(p-1)w_n=\min_{\sigma\in S_n}\sum_{i=1}^{n} \pceil{\frac{pi-\sigma(i)}{d}}.
\end{equation*}
Write $u_n=v_1\cdots v_n$, which does not depend on $\sigma$,  
so we have
\begin{align*}
\det((\bG_{pi-j})_{1\le i,j\le n})
&=\sum_{\sigma\in S_n}\sgn(\sigma)\prod_{i=1}^n \bG_{pi-\sigma(i)}\\
&\equiv
u_n\sum_{\sigma\in S_n}
\sgn(\sigma)A_d^{\sum_{i=1}^n \pceil{\frac{pi-\sigma(i)}{d}}-nd}
\prod_{i=1}^{n}\tH_{i,\sigma(i)}\; \gamma^{\sum_{i=1}^n \pceil{\frac{pi-\sigma(i)}{d}}}
  \\
& \stackrel{(\star)}{\equiv} u_n\sum_{\sigma\in S_n^{\min}}
\sgn(\sigma)A_d^{N_n}
\prod_{i=1}^{n}\tH_{i,\sigma(i)} \gamma^{(p-1)w_n}
\\
&\equiv 
u_n\gamma^{(p-1)w_n}A_d^{N_n}\widetilde{f}_n 
 \pmod {\gamma^{(p-1)w_n+1}},
\end{align*}
where ($\star$) is due to the observation that 
$\sum_{i=1}^n \pceil{\frac{pi-\sigma(i)}{d}} = (p-1)w_n$ if and only if 
$\sigma\in S_n^{\min}$.
\end{proof}

\begin{prop}
\label{P:est}
Let $p>(d^2+1)(d-1)$.
If $\vc=(c_1,\ldots,c_d)\in {\bar\Z_p}^d$ such that $v_p(P_r(\vc))=0$ for each factor $P_r$ of 
$\tilde\Psi_r$,
 then we have
$$\left\{
\begin{array}{lll}
v_p(\bG_{pi-j}(\vc))&=&\frac{\pceil{\frac{pi-j}{d}}}{p-1},\\
v_p\left(\det(\bG_{pi-j}(\vc))_{i,j\le n}\right) &=&w_n. 
\end{array} 
\right.
$$
\end{prop}

\begin{proof}
By the hypothesis, $v_p(c_d)=v_p(\tH_{ij}(\vc))=0$, so we have by the first congruence of 
Lemma \ref{L:congruence}
\begin{align*}
v_p(\bG_{pi-j}(\vc))
=
v_p\left(v_i\gamma^{\pceil{\frac{pi-j}{d}}}c_d^{\pceil{\frac{pi-j}{d}-d}}\tH_{ij}(\vc)\right)
=\pceil{\frac{pi-j}{d}}v_p(\gamma)=\frac{\pceil{\frac{pi-j}{d}}}{p-1}.
\end{align*}
This proves the first valuation. 
The second one follows the same argument
using also the hypothesis $v_p(\tf_n(\vc))=0$ 
and the second congruence of Lemma \ref{L:congruence},
hence is omitted.
\end{proof}

\section{Main theorems and their applications}
\label{S:proofs}

This section proves our main results. 
By \eqref{E:zeta_function}, 
$X_{f\bmod\wp,\ell}$ is ordinary (see Definition \ref{D:ord}) if and only if 
$\NP_q(L_{f\bmod \wp}(\chi_\ell,s))=\GNP(\A^d,\F_p)$.
We shall first prove Theorem \ref{T:main2} that concerns $L$-function 
$L_{f\bmod \wp}(\chi_\ell,s)$. 
Then we apply this result to the Zeta function and prove Theorem \ref{T:main1b}.

We now have to introduce essential notations for our theorems before we proceed.
Define 
\begin{eqnarray}
\label{E:cU}
\cU &\coloneqq &\bigcap_{r\in(\Z/d\Z)^*}\tilde{\cW_r},
\end{eqnarray}
where $\tilde\cW_r$ is from Proposition \ref{P:open}.
As we have proved there, 
it is Zariski dense open in $\A^d$, defined over $\Q$, and depends only on $d$.

Starting in this section, we will need to work with $\bar\Q$ and $\bar\Q_p$ 
(for a prime $p$) simultaneously, so we  
shall tacitly fix an embedding from $\bar\Q$ to $\bar\Q_p$ once and for all.  

\begin{definition*}
Suppose $f=\sum_{m=1}^d a_mx^m\in\cU(\bar\Q)$ and denote $\va=(a_1,\ldots,a_d)$.
Let $\heit(f)$ be the least integer $N\ge (d^2+1)(d-1)$
such that for every prime $p_0>N$ we have the following properties satisfied
\begin{enumerate}
\item[i)] $f$ is $p_0$-adically integral, i.e.,
$v_{p_0}(a_m)\ge 0$ for $m=1,\ldots,d$;
\item[ii)] 
If $p_0\equiv r\bmod d$, then
$P_r(\va)$ is a $p_0$-adically unit, i.e., 
$v_{p_0}(P_r(\va))= 0$, for every factor $P_r$ of $\tilde\Psi_r$ (see Definition \ref{D:Psi_r}).
\end{enumerate}
We call $\heit(f)$ the {\bf height of $f$}.
\end{definition*}
First of all let us note that $\heit(f)$ is well-defined: For any prime $p_0>d$, by Proposition \ref{P:open}(2), $P_r\in\Z_{p_0}[A_1,\ldots,A_d]$ is nonzero. 
When $p_0$ is large enough so that 
it is coprime to the denominator of each $a_m\in \bar\Q$, then  
we have i) satisfied. Since $f\in\cU(\bar\Q)$, 
$P_r(\va)\in \bar\Q^*$, hence property ii) is satisfied for $p_0$ large enough. 
Secondly, we remark that $\heit(f)$ depends only on $f$.

For any prime $\wp$ of $K$, we write $k(\wp)$ for
its residue field. 
When $\va\in {\bar\Z_p}^d$ and $\wp|p$, 
we denote $\bar{a}_m=(a_m\bmod \wp)$ in $k(\wp)$, and 
write $\hat\va=(\hat{a}_1,\ldots,\hat{a}_d)$ for its Teichm\"uller lifting.

\begin{lemma}\label{L:lemma1}
Suppose  $f\in\cU(K)$ for a number field $K$ and suppose $r\in(\Z/d\Z)^*$.
If $p\equiv r\bmod d$ and $p>\heit(f)$, then
we have 
$$v_p(P_r(\hat\va))=0$$ 
for every factor $P_r$ of $\tilde\Psi_r$.
\end{lemma}

\begin{proof}
By Proposition \ref{P:open}(2), $P_r\in\Z_p[A_1,\ldots,A_d]$ is nonzero since $p>\heit(f)>d$.
Since $\hat{\va}\in(\bar\Z_p)^d$,  $v_p(P_r(\hat\va))\ge 0$. 
Let $\wp|p$. Since $\hat\va\equiv \va \bmod \wp$, we have 
$P_r(\hat\va)\equiv P_r(\va)\bmod \wp$. 
Therefore, $v_p(P_r(\hat\va))=0$ if and only if  
$v_p (P_r(\va))=0$. 
On the other hand, 
since $p>\heit(f)$, $v_p(P_r(\va))=0$. 
Thus $v_p(P_r(\hat\va))=0$. 
\end{proof}

\subsection{Asymptotic $L$-function lower bound}

For every $\bar{f}\in\A^d(\F_q)$ we have
$\NP_q(L_{\bar\alpha \bar{f}}(\chi_1,s))\ge \GNP(\A^d,\F_p)$
for all $\bar\alpha\in\F_q^*$ (see \cite{Zh04}). 
The equality can be achieved for many 
polynomials $\bar{f}\in \A^d(\F_q)$, 
but this may depend on the prime $p$ and $\bar\alpha$.

\begin{theorem}
\label{T:sameNP}
Let $f\in\cU(K)$ for a number field $K$. 
For every prime $\wp$ of $K$
such that $p=\chara(k(\wp))>\heit(f)$, 
we have 
$$ 
\NP(L_{\bar\alpha (f\bmod \wp)}(\chi_1,s))=\GNP(\A^d,\F_p)
$$
for every $\bar\alpha\in k(\wp)^*$.
\end{theorem}

\begin{proof}
Write $f=\sum_{m=1}^{d} a_m x^m$, 
then all $a_m$ are $p$-adic integral since $p>\heit(f)$. 
Write $q=p^b=|k(\wp)|$.
Then $\bar{f}=(f\bmod\wp)$ is over $\F_q$.
To simplify notation, for the rest of this proof (only) 
let $G=(G_{ij})_{i,j\ge 1}$ where
$G_{ij}:=\bG_{pi-j}(\hat\alpha\hat\va)$ lies in $\Q_q(\zeta_p)$.
Then 
\begin{eqnarray}\label{E:GGG}
G &=& G(\bar\alpha\bar{f})=(\bG_{pi-j}(\hat\alpha\hat\va))_{i,j\ge 1}
\end{eqnarray}
in the notation of \eqref{E:Gfbar}. 

First of all we shall verify  all hypotheses of Proposition \ref{P:compute!} for 
$G=(G_{ij})_{i,j\ge 1}$ are satisfied.
Since $d\ge 3$, we have $p>(d^2+1)(d-1)>2d^2-2d+1$.
First of all, by \eqref{E:G_n}, 
$G_{ij}=\bG_{pi-j}(\hat\alpha\hat\va)=\sum(\prod_{i=1}^d\lambda_{m_i}){\hat\alpha}^{\sum_{i=1}^d m_i}\prod_{i=1}^d\hat{a}_i^{m_i}$ with the sum ranges over 
all nonnegative integers $m_1,\ldots,m_d$ such that $\sum_{i=1}^d i m_i =pi-j$. 
Since $v_p(\lambda_{m_i})\ge \frac{m_i}{p-1}$, $v_p(\hat\alpha)=v_p(\hat{a}_d)=0$, 
and $v_p(\hat{a}_i)=0$ for all $a_i\ne 0$,
we have $v_p(G_{ij})\ge 
\frac{\pceil{\frac{pi-j}{d}}}{p-1}\ge \frac{pi-j}{d(p-1)}$ for all $i,j\ge 1$.
This proves property (1) holds.
As $f\in\cU(K)$ and $p>\heit(f)$, 
by Lemma \ref{L:lemma1}, we have $v_p(P_r(\hat\va))=0$ for every factor $P_r$ of
$\tilde\Psi_r$. Since $P_r$ is homogenous (denote its total degree $d_r$) 
we have $P_r(\hat\alpha \hat\va)=\hat\alpha^{d_r}P_r(\hat\va)$.
Thus $$v_p(P_r(\hat\alpha \hat\va))=v_p(\hat\alpha^{d_r})+v_p(P_r(\hat\va))
=v_p(P_r(\hat\va))=0.$$ This allows us to apply Proposition \ref{P:est} and conclude that
properties (2) and (3) are satisfied.
With all the above, we now may apply Proposition \ref{P:compute!} to the matrix 
$G$ to get 
\begin{eqnarray}\label{E:GGG2}
\NP_q^{<1}(\det(1-G_{[b]}s))=\GNP(\A^d,\F_p).
\end{eqnarray}
Secondly, since \eqref{E:GGG},
we apply Proposition \ref{P:L-function} to $\bar\alpha \bar{f}$ and have 
\begin{eqnarray}\label{E:GGG3}
\NP_q(L_{\bar\alpha\bar{f}}(\chi_1,s))= \NP_q^{<1}(\det(1-G(\bar{\alpha}\bar{f})_{[b]}s))=
 \NP_q^{<1}(\det(1-G_{[b]}s)).
\end{eqnarray}
Combining \eqref{E:GGG2} and \eqref{E:GGG3}, we get 
$\NP_q(L_{\bar\alpha(f\bmod \wp)}(\chi_1,s))=\GNP(\A^d,\F_p)$.
\end{proof}

\begin{definition}
\label{D:HP}
The {\bf Hodge polygon} $\HP(\A^d)$ for $\A^d$ is the Newton polygon  
with vertices 
at $$(0,0), \left(1,\frac{1}{d}\right), \cdots, \left(n,\frac{n(n+1)}{2d}\right), \ldots,\left(d-1,\frac{(d-1)d}{2d}\right).$$
 Notice that $\GNP(\A^d,\F_p)\ge \HP(\A^d)$ and 
they are equal only when $p\equiv 1\bmod d$.
The {\bf Hodge polygon} for genus-$g$ Artin-Schreier curves 
is the dilation of $\HP(\A^d)$ by a factor of $p^\ell-1$, we denote it by $\HP_g$.  
\end{definition}

It is well known that $\HP(\A^d)$ is the lower bound 
for all $\NP_q(L_{\bar{f}}(\chi_\ell,s))$ for all $\bar{f}\in\A^d(\bar\F_p)$ (see \cite{AS87}\cite{Wan93}\cite{Zh04}). For any $\bar{f}\in\A^d(\bar\F_p)$ 
we always have 
\begin{equation*}
\NP_q(L_{\bar{f}}(\chi_\ell,s))\hh\ge \hh\GNP(\A^d,\F_p)\hh\ge \hh \HP(\A^d).
\end{equation*}
The following theorem generalizes the rank-1 case (that is, $\ell=1$) proved in \cite{Zh04}.

\begin{theorem}
\label{T:main2} \-
\begin{enumerate}
\item For every $f\in\A^d(K)$ and any prime $\wp|p$ of a number field $K$.
Then we have
\begin{equation*}
\NP_q(L_{f\bmod \wp}(\chi_\ell,s))
\hh\ge \hh \GNP(\A^d,\F_p)\hh \ge \hh \HP(\A^d).
\end{equation*}
\item  If $f\in\cU(K)$,
then 
\begin{equation*}
\NP_q(L_{f\bmod \wp}(\chi_\ell,s))
\hh = \hh \GNP(\A^d,\F_p)
\end{equation*}
for every prime $\wp$ of $K$ with 
$p=\chara(k(\wp))>\heit(f)$ and $\deg(\wp)\in\ell\Z$.
\end{enumerate}
\end{theorem}

\begin{proof}
(1) By Proposition \ref{P:reduction}, 
\begin{equation}\label{E:L_L}
L_{f\bmod \wp}(\chi_\ell,s) \hh = \hh L_{\bar{\alpha}(f\bmod \wp)}(\chi_1,s)
\end{equation}
for some  $\bar{\alpha}\in\F_{p^\ell}^*$. But we always have (see for example \cite{Zh04}) 
\begin{equation*}
\NP_q(L_{\bar{\alpha}(f\bmod \wp)}(\chi_1,s)) \hh\ge \hh \GNP(\A^d,\F_p).
\end{equation*}
This proves the first statement:
\begin{equation*}
\NP_q(L_{f\bmod \wp}(\chi_\ell,s)) \hh \ge \hh\GNP(\A^d,\F_p).
\end{equation*}

(2) Now let $f\in\cU(K)$.
Note that 
$\bar\alpha\in\F_{p^\ell}\subseteq \F_{p^{\deg(\wp)}}=k(\wp)$.
By Theorem \ref{T:sameNP}, we have 
$
\NP_{q}(L_{\bar\alpha(f\bmod \wp)}(\chi_1,s))=\GNP(\A^d,\F_p).
$
Combining this with above \eqref{E:L_L}, we have
$$\NP_{q}(L_{(f\bmod \wp)}(\chi_\ell,s))
=\NP_{q}(L_{\bar\alpha(f\bmod \wp)}(\chi_1,s))
=\GNP(\A^d,\F_p).$$
This concludes the proof.
\end{proof}

\begin{corollary}
\label{C:indep}
If $f\in \cU(K)$ for a number field $K$
and if $\wp$ is any prime of $K$ with $p=\chara(k(\wp))>\heit(f)$ and $\deg(\wp)\in\ell\Z$,
then the normalized Newton polygon of $L_{f\bmod\wp}(\chi_\ell,s)$ 
is independent of the choice of 
the irreducible, non-trivial character $\chi_\ell:\F_{p^\ell}\to \bar\Q_p^*$.
\end{corollary}

\begin{proof}
By Theorem \ref{T:main2},  
the normalized Newton polygon 
of $L$-function $L_{f\bmod \wp}(\chi_\ell,s)$ is always equal to $\GNP(\A^d,\F_p)$, independent of the choice of $\chi_\ell$.
\end{proof}

\begin{corollary}
\label{C:limit}
Let $f\in\cU(K)$. Let $p>\heit(f)$.  
If $\wp|p$ is a prime of $K$ with $\deg(\wp)\in\ell\Z$,
then we have 
$$
\lim_{p\rightarrow \infty} 
\NP_{|k(\wp)|}
(L_{(f\bmod \wp)}({\chi_\ell},s)) 
= \HP(\A^d).
$$
\end{corollary}

\begin{proof}
By Theorem \ref{T:main2}, for $p$ large enough the  Newton polygon 
of $L_{(f\bmod \wp)}(\chi_\ell,s)$ is equal to $\GNP(\A^d,\F_p)$.
Notice that $\lim\limits_{p\ra \infty}w_n=\frac{n(n+1)}{2d}$ 
for every $1\le n\le d-1$. 
Hence we have $\lim\limits_{p\ra \infty}\GNP(\A^d,\F_p)=\HP(\A^d)$, 
and this proves our statement. 
\end{proof}

 \subsection{Asymptotic lower bound for Zeta functions of curves}
\label{S:theorems}
 
It is already known (see \cite{Zh04} for references)
that for $p$ large enough the lower bound of all $\NP(L_{\bar{f}}(\chi_1,s))$ can be achieved and is $\GNP(\A^d,\F_p)$ (as given in Definition \ref{D:GNP}).
We shall show below this generalizes to rank-$\ell$ case.
Recall $\GNP_{g,\ell,\bar\F_p}$ from Definition \ref{D:ord}.

\begin{theorem}[Theorem \ref{T:main1}]
\label{T:main1b}
Let $K$ be a number field.
\begin{enumerate}
\item
For every $f\in\A^d(K)$ and every prime $\wp|p$
with $\deg(\wp)\in\ell\Z$, 
we have 

$$
\NP(X_{(f\bmod \wp),\ell})\ge \GNP(\A^d,\F_p)^{p^\ell-1}\ge \HP_g.
$$
\item
If  $f\in \cU(K)$, 
then every prime $\wp$ of $K$ with $\chara(k(\wp))
>\heit(f)$ and $\deg(\wp)\in\ell\Z$ is ordinary, that is, 
    \begin{equation*}
        \NP(X_{(f\bmod \wp),\ell}) = \GNP(\A^d,\F_p)^{p^\ell-1}=\GNP_{g,\ell,\bar\F_p}.
    \end{equation*}
\end{enumerate}
\end{theorem}

\begin{proof}
(1). 
Write $q=|k(\wp)|$. 
Applying Theorem \ref{T:main2}(1) to $\bar\alpha(f\bmod \wp)$,   we get 
$\NP_q(L_{\bar\alpha(f\bmod \wp)}(\chi_1,s))\ge \GNP(\A^d,\F_p)$.
By our hypothesis that $\deg(\wp)\in\ell\Z$
we have $\F_{p^\ell}\subseteq \F_{p^{\deg(\wp)}}=k(\wp)$. 
So we may apply \eqref{E:zeta_function} and get
$$\NP(X_{f\bmod \wp,\ell})=\bigcup_{\bar\alpha\in\F_{p^\ell}^*}\NP_q(L_{\bar\alpha(f\bmod \wp)}(\chi_1,s))\ge \GNP(\A^d,\F_p)^{p^\ell-1},
$$
hence we have the first inequality.
Since $\GNP(\A^d,\F_p)\ge \HP(\A^d)$, we have $$\GNP(\A^d,\F_p)^{p^\ell-1}\ge \HP(\A^d)^{p^\ell-1}=\HP_g,$$ this proves the second inequality.

\noindent 
(2). 
Applying Theorem \ref{T:sameNP}, for every $f\in\cU(K)$ and every $\bar\alpha\in\F_{p^\ell}^*$
we have $$\NP_q(L_{\bar\alpha(f\bmod \wp)}
(\chi_1,s))=\GNP(\A^d,\F_p).$$
Combined with \eqref{E:zeta_function}, we have 
$$
\NP(X_{(f\bmod \wp),\ell})
= \bigcup_{\bar\alpha\in\F_{p^\ell}^*}\NP_q(L_{\bar\alpha(f\bmod \wp)}(\chi_1,s))
=\GNP(\A^d,\F_p)^{p^\ell-1}.$$
Thus we have 
$\GNP(\A^d,\F_p)^{p^\ell-1}=\GNP_{g,\ell,\bar\F_p}$.
This proves the theorem.
\end{proof}

\begin{corollary}[Corollary \ref{C:1}] 
\-
\begin{enumerate}
\item  Let $f\in\cU(K)$ for any number field $K$. 
If $\wp$ is a prime of $K$ with $\deg(\wp)\in\ell\Z$, and 
$p=\chara (k(\wp))>\heit(f)$ and $r=(p\bmod d)$,
then the slope multi-set of
$\NP(X_{f\bmod \wp,\ell})$ is equal to 
$$
\left\{\frac{1}{d}+\frac{\varepsilon_1}{d(p-1)},\frac{2}{d}+\frac{\varepsilon_2}{d(p-1)},\ldots, \frac{d-1}{d}+\frac{\varepsilon_{d-1}}{d(p-1)}\right\}^{p^\ell-1},$$
where $\varepsilon_i\in \Z$ depends only on $r$ and $d$, and $-(i-1)(d-1)\le \varepsilon_i \le i(d-1)$.

\item  When $p$ is large enough, $\GNP_{g,\ell,\bar\F_p}=\HP_g$ if and only if $p\equiv 1\bmod d$. 
\end{enumerate}
\end{corollary}

\begin{proof}
(1) By Definition \ref{D:GNP}, the Newton polygon $\GNP(\A^d,\F_p)$ 
has exactly $d-1$ line segments. 
The $i$-th segment is of slope
$
w_i-w_{i-1}
=\frac{i}{d}+\frac{\varepsilon_i}{d(p-1)}
$
where $\varepsilon_i=M_i-M_{i-1}$.
It is clear from \eqref{E:M_n} that 
$0\le M_i\le i(d-1)$, thus 
$-(i-1)(d-1)\le \varepsilon_i
\le i(d-1)$.
Then $\NP(X_{(f\bmod \wp),\ell})$ is precisely $p^\ell-1$ copies of the above slopes,
and our statement follows from
Theorem \ref{T:main1} above.\\
(2) 
For $p$ large enough then we have $\GNP_{g,\ell,\bar\F_p}$ is equal to 
$\NP(X_{f\bmod \wp,\ell})$ for some $f\in \cU(\bar\Q)$, and its slopes are demonstrated in part (1). It remains to show that $p\equiv 1\bmod d$ if and only if $M_i=M_{i-1}$ for all $i=1,\ldots,d$ where we set $M_0=0$. It is a straightforward computation which we leave out as an exercise to the interested reader.
\end{proof}

\begin{corollary}\label{C:2}
Let $K$ be a number field that contains a $\Z/\ell\Z$-subextension of $\Q$.
Let $f\in\cU(K)$,
then $X_{f\bmod\wp,\ell}$ is ordinary 
for infinitely many primes $\wp$ of $K$. 
\end{corollary}
\begin{proof}
By hypothesis on $K$ and class field theory, 
there are primes $\wp|p$ with $\deg(\wp)\in\ell\Z$ for a congruence class of $p$ 
(hence infinitely many). 
Now apply Theorem \ref{T:main1b} above, 
we find that all primes $\wp|p$ with $p>\heit(f)$ are ordinary. 
Hence there are infinitely many such $\wp$.
\end{proof}

\section{1-parameter family of $f(x)=x^d+ax$}
\label{S:6}

Let $\A^d(1)$ denote the coefficient space of all $f=x^d+ax$. 
This is a 1-parameter subspace of $\A^d$.
It is known (see \cite[Section 6]{Zh03}) 
that for any $f=x^d+ax$ with $a\neq 0$, we have
$\NP(L_{f\bmod \wp}(\chi_1,s))=\GNP(\A^d(1),\F_p)$
for all prime $\wp|p$ where $p$ is large enough.
Even though $\GNP(\A^d(1),\F_p)\ge \GNP(\A^d,\F_p)$ and they are not 
equal in general, 
they always have the same limit as 
$p\rightarrow\infty$:
$$
\lim_{p\rightarrow \infty}\GNP(\A^d(1),\F_p)=\lim_{p\rightarrow\infty}\GNP(\A^d,\F_p)=\HP(\A^d).
$$
In this section we prove that this result above generalizes to 
rank-$\ell$ $L$-function. 

\begin{lemma}\label{L:1-para}
Let $p> (d-1)^3+1$ and $p\equiv r\bmod d$. For every $1\le i,j,n\le d-1$ we have 
$$\left\{
\begin{array}{lcll}
\bG_{pi-j}
&\equiv &v'_{i,j}A_1^{r'_{i,j}}A_d^{\pfloor{\frac{pi-j}{d}}}\gamma^{\pfloor{\frac{pi-j}{d}}+r'_{ij}}& (\bmod \gamma^{\pfloor{\frac{pi-j}{d}}+r'_{ij}+1});\\
\det((\bG_{pi-j})_{1\le i,j\le n})
&\equiv & u'_nA_1^{M'_n}A_d^{\frac{(p-1)n(n+1)}{2d}-\frac{M'_n}{d}}\gamma^{(p-1)w'_n} 
& (\bmod \gamma^{(p-1)w'_n+1})
\end{array}
\right.
$$
for some $v'_{i,j}$ and $u'_n\in\Z_p^*$.
\end{lemma}

\begin{proof}
Computation in this proof is similar to that in Proposition \ref{P:est}.
Write $v'_{i,j}=(\pfloor{\frac{pi-j}{d}}! r'_{i,j}!)^{-1}$, it is easy to see 
that for $p>d$ the first congruence holds.
Let 
$u'_n=\det((\pfloor{\frac{pi-j}{d}}!r'_{ij}!)^{-1})$.
For $p>(d-1)^3+1$ the second congruence holds. To show that $u'_n$ is a $p$-adic unit, we observe that
$$u_n'=\prod_{i=1}^n (\pfloor{\frac{pi-1}{d}}! r'_{i1}!)^{-1}M'$$
where $M'=\det(r'_{k1}(r'_{k1}-1)\cdots (r'_{k1}-m+1))_{1\le k,m\le n}$.
Notice that $M'=\prod_{1\le k<m\le n} (r'_{k1-r'_{m1}})\in\Z_p^*$ and 
hence $u_n'\in\Z_p^*$. This proves our second congruence.
\end{proof}

For any $a\in\bar\Q$, let $\denom(a)$ be the largest prime 
factor in the denominator of $N_{\Q(a)/\Q}(a)$. 
The following theorem is an analog of Theorem \ref{T:main2}.

\begin{theorem}
\label{T:1-parameter}
Let $f=x^d+ax\in\A^d(\bar\Q)$ with $a\ne 0$.
\begin{enumerate}
\item For any prime $\wp|p$ with $p>\max((d-1)^3+1, \denom(a))$ 
and $\deg(\wp)\in\ell\Z$, 
we have 
\begin{eqnarray*}
\NP(L_{\bar\alpha(f\bmod \wp)}) 
(\chi_1,s))&=&\GNP(\A^d(1),\F_p) \quad {\mbox{for any $\bar\alpha\in\F^*_{p^\ell}$}},
\\
\NP(L_{f\bmod \wp}(\chi_\ell,s))&=&\GNP(\A^d(1),\F_p).
\end{eqnarray*}
\item  For every prime $p$, let $\wp|p$ be a prime in a number field $K\supseteq \Q(a)$ 
with $\deg(\wp)\in\ell\Z$, then  
$$\lim_{p\rightarrow\infty} \NP(L_{f\bmod \wp}(\chi_\ell,s))=
\lim_{p\rightarrow \infty}\GNP(\A^d(1),\F_p) 
=\HP(\A^d).$$
\end{enumerate}
\end{theorem}

\begin{proof}
(1) By the same argument as that of 
Theorem \ref{T:main2}(2), we notice that the second equality follows from the first. 
The first equality can be proved by following the same argument of 
Theorem \ref{T:sameNP}.
Let $G_{ij}:=\bG_{pi-j}(\hat\alpha\hat{a})$.
All is reduced to show that the matrix $G=(G_{ij})_{i,j\ge 1}$ 
satisfies the hypothesis of Proposition \ref{P:compute!} (the second case). 
Indeed, we see easily that
$v_p(G_{ij})\ge \frac{\pfloor{\frac{pi-j}{d}}+r_{ij}'}{p-1}
\ge\frac{pi-j}{d(p-1)}$ for all $i,j\ge 1$.
Then by Lemma \ref{L:1-para}
for $p>(d-1)^3+1$ 
and $1\le i,j,n\le d-1$,
we have for $a\ne 0$ that 
$$\left\{
\begin{array}{lcl}
v_p(G_{ij}) &=& \frac{\pfloor{\frac{pi-j}{d}}+r_{ij}'}{p-1},\\
v_p(\det(G_{ij})_{1\le i,j\le n})&=&w'_n.
\end{array}
\right.
$$
By applying Proposition \ref{P:compute!} to the matrix 
$G=(G_{ij})_{i,j\ge 1}$ with $h_i=\frac{i}{d}$ and 
$\delta=\frac{1}{d}$, this verifies the hypotheses are satisfied and we have
proved the first part of the theorem.

(2) The second part is straightforward as 
$\lim_{p\ra \infty}w_n'=\frac{n(n+1)}{2d}$.
\end{proof}

To conclude, we have an analog theorem of Theorem \ref{T:main1} (or Theorem \ref{T:main1b}) 
for the family $\A^d(1)$, that is, for every $f=x^d+ax$ and $a$ is nonzero in $\bar\Q$,
the Artin-Schreier curve $X_{f\bmod \wp,\ell}$ 
is ordinary for all $p$ large enough and $\wp|p$ with $\deg(\wp)$ divisible by $\ell$. 

\begin{theorem}
\label{T:main1c}
Let $K$ a number field with $a\in K^*$.
If $f=x^d+ax$, then we have $$\NP(X_{f\bmod\wp,\ell})=\GNP(\A^d(1),\F_p)^{p^\ell-1}$$  
for all primes $\wp|p$ with $\deg(\wp)\in\ell\Z$ and $p>\max((d-1)^3,den(a))$. 
\end{theorem}

\begin{proof}
This is an immediate consequence from Theorem \ref{T:1-parameter},
following the same argument as that of Theorem \ref{T:main1b}. 
\end{proof}

\end{document}